# CORRECTION

# ERROR ESTIMATES FOR BINOMIAL APPROXIMATIONS OF GAME OPTIONS

BY YURI KIFER

My student Y. Dolinsky noticed that the inequalities (5.33) needed to obtain the estimate (5.34) in the proof of Theorem 2.3 hold true only for hedging strategies without short selling of bonds and stocks, that is, when the amounts of bonds and stocks in the portfolio are always nonnegative. Since the existence of such hedging strategies cannot be guaranteed, in general, if their initial capital equals the option price, the proof should be corrected and we start with an argument due to Dolinsky which serves this purpose. In the notation of [1] set

$$\Psi = \sup_{0 \le t \le T} \left( Q_z^B(\theta_\varphi^{(n)}, t) - Q_z^{B,n}\left(\frac{\varphi T}{n}, \frac{\nu_t T}{n}\right) \right)^+. \tag{1}$$

From (5.29)–(5.32) of [1], we obtain that there exists a constant $C > 0$ such that

$$E^B \Psi \le C(F_0(z) + \Delta_0(z) + z + 1)n^{-1/4}. \tag{2}$$

Let $\tau \in \mathcal{T}_{0T}^B$ be a stopping time. Then $\nu_\tau = \min\{k \in \mathbb{N} : \theta_k^{(n)} \ge \tau\} \in \mathcal{T}^{B,n}$, and so $\theta_{\nu_\tau}^{(n)} \in \mathcal{T}^B$ is a stopping time (see beginning of proofs of Lemmas 3.1 and 3.6 in [1]). As any self-financing discounted portfolio $\check{Z}^B$ (see (5.21) in [1]) is a martingale, and so taking into account (5.24) in [1] and the optional sampling theorem, we derive that

$$\check{Z}_{\theta_\varphi^{(n)} \wedge \tau}^B = E^B(\check{Z}_{\theta_\varphi^{(n)} \wedge \theta_{\nu_\tau}^{(n)}}^B | \mathcal{F}_{\theta_\varphi^{(n)} \wedge \tau})$$

$$\ge E^B\left( Q_z^{B,n}\left(\frac{\varphi T}{n}, \frac{\nu_\tau T}{n}\right) \Big| \mathcal{F}_{\theta_\varphi^{(n)} \wedge \tau} \right)$$
(3)

---

Received October 2006.

$$\geq E^B(Q_z^B(\theta_\varphi^{(n)},\tau) - \Psi|\mathcal{F}_{\theta_\varphi^{(n)}\wedge\tau})$$
$$= Q_z^B(\theta_\varphi^{(n)},\tau) - E^B(\Psi|\mathcal{F}_{\theta_\varphi^{(n)}\wedge\tau}).$$

Finally, (3) yields that

$$\sup_{\tau\in\mathcal{T}_{0T}^B} E^B(Q_z^B(\theta_\varphi^{(n)},\tau) - \check{Z}_{\theta_\varphi^{(n)}\wedge\tau}^B)^+ \leq E^B\Psi,$$

which together with (2) provides the required estimate for the expectation of the left-hand side of (5.26) in [1] and (2.22) there follows. □

There is another way to fix (5.33) in [1] which is interesting by itself. Since the main point of Theorem 2.3 in [1] is to show how to construct in an explicit way some "nearly" hedging strategies for the Black–Scholes market using hedging strategies in approximating Cox, Ross, Rubinstein (CRR) markets, it suffices, essentially, to consider only hedging strategies in CRR markets which are given by standard explicit formulas via the Doob decomposition and the discrete martingale representation in CRR models. Thus, we modify (5.33) in [1] by writing first (in notation of [1]),

$$(4)\quad (\check{Z}_{\theta_\varphi^{(n)}\wedge\theta_{\nu_t}^{(n)}}^B - \check{Z}_{\theta_\varphi^{(n)}\wedge t}^B)^+ \mathbb{I}_{A_t} \leq \max_{1\leq k\leq n}|\gamma_{\theta_k^{(n)}}^\varphi|(\check{S}_{\theta_\varphi^{(n)}\wedge\theta_{\nu_t}^{(n)}}^B(z) - \check{S}_{\theta_\varphi^{(n)}\wedge t}^B(z))^+$$

and then estimating $\gamma_{\theta_k^{(n)}}^\varphi$ via explicit formulas. Namely (see, e.g., (2.28) in Theorem 1 of [3]),

$$(5)\qquad \gamma_{\theta_k}^\varphi = \frac{\alpha_k e^{rTk/n}}{S_{(k-1)T/n}^{B,n}(z)},$$

where $S^{B,n}$ is defined in (3.1) of [1] and $\alpha_k$ comes from the explicit martingale representation formula (see Section 3 in [3] or Section 4d, Chapter V in [2])

$$(6)\qquad M_k = M_0 + \sum_{j=1}^k \alpha_j(\rho_j^{(n)} - r^{(n)}).$$

Here $\rho_j^{(n)}$ and $r^{(n)}$ are given by (1.8) in [1] as

$$r^{(n)} = e^{rT/n} - 1 \quad\text{and}\quad \rho_j^{(n)} = \exp(rT/n + \kappa(B_{\theta_j^{(n)}}^* - B_{\theta_{j-1}^{(n)}}^*)), \qquad j\geq 1,$$

and $\{M_k\}_{0\leq k\leq n}$ is the martingale emerging from the Doob decomposition of the supermartingale

$$(7)\quad U_k = U_k^\zeta = \max_{\nu\in\mathcal{J}_{0,n}} E(Q_z^{B,n}(\zeta T/n, (\nu\circ\lambda_B^{(n)})T/n)|\mathcal{G}_k^{B,n}), \qquad \zeta\in\mathcal{S}_{0,n}^{B,n},$$

where the notation is the same as in [1]. The martingale in (3) can be written explicitly in the form (see Section 1b, Chapter II in [2])

$$(8) \qquad M_k = U_0 + \sum_{j=1}^{k}(U_j - E(U_j|\mathcal{G}_{j-1}^{B,n})).$$

Since $U_k$ in (7) is measurable with respect to

$$\mathcal{G}_k^{B,n} = \sigma(B^*_{\theta_1^{(n)}}, B^*_{\theta_2^{(n)}} - B^*_{\theta_1^{(n)}}, \ldots, B^*_{\theta_k^{(n)}} - B^*_{\theta_{k-1}^{(n)}}),$$

we can write also

$$(9) \qquad U_k = \Psi_k(B^*_{\theta_1^{(n)}}, B^*_{\theta_2^{(n)}} - B^*_{\theta_1^{(n)}}, \ldots, B^*_{\theta_k^{(n)}} - B^*_{\theta_{k-1}^{(n)}}),$$

where, in view of the assumption (2.1) from [1], the functions $\Psi_k$ satisfy

$$(10) \qquad |\Psi_k(x_1,\ldots,x_k) - \Psi_k(y_1,\ldots,y_k)|$$
$$\leq Cz \max_{0\leq j\leq k}\left|\exp\left(\kappa\sum_{i=1}^{j}x_i\right) - \exp\left(\kappa\sum_{i=1}^{j}y_i\right)\right|$$

for some $C = C_T > 0$ independent of $k, n, \{x_j\}$ and $\{y_j\}$ but depending on $T$. Observe that $\alpha_k$ can also be written in the explicit form (see Section 3 in [3] or Section 4d, Chapter V in [2]) which in our situation amounts to

$$(11) \qquad \begin{aligned}\alpha_k &= e^{-rT/n}(e^{\kappa\sqrt{T/n}} - 1)^{-1}(1 - p_n)\\ &\quad \times (\Psi_k(B^*_{\theta_1^{(n)}}, B^*_{\theta_2^{(n)}} - B^*_{\theta_1^{(n)}}, \ldots, B^*_{\theta_{k-1}^{(n)}} - B^*_{\theta_{k-2}^{(n)}}, \sqrt{T/n})\\ &\quad - \Psi_k(B^*_{\theta_1^{(n)}}, B^*_{\theta_2^{(n)}} - B^*_{\theta_1^{(n)}}, \ldots, B^*_{\theta_{k-1}^{(n)}} - B^*_{\theta_{k-2}^{(n)}}, -\sqrt{T/n})),\end{aligned}$$

where $p_n = (e^{\kappa\sqrt{T/n}} - 1)^{-1}$. This together with (7) yields that

$$(12) \qquad |\alpha_k| \leq 2Cz\exp(\kappa B^*_{\theta_{k-1}^{(n)}}).$$

Now by (2) and (3.1) of [1],

$$(13) \qquad \max_{1\leq k\leq n}|\gamma^\varphi_{\theta_k^{(n)}}| \leq 2C.$$

Finally, estimating $(\check{S}^B_{\theta_\varphi^{(n)}\wedge\theta_{\nu_t}^{(n)}}(z) - \check{S}^B_{\theta_\varphi^{(n)}\wedge t}(z))^+$ in the same way as in the last line of (5.33) in [1] and using (13), we arrive at (5.34) in [1] in the same way as there and Theorem 2.3 follows. $\square$

We warn the reader also about the misprint in (4.42) of [1] where $t$ before $Q_z^{B,n,\theta}$ should be deleted.

## REFERENCES

[1] KIFER, YU. (2006). Error estimates for binomial approximations of game options. *Ann. Appl. Probab.* **16** 984–1033. MR2244439
[2] SHIRYAEV, A. N. (1999). *Essentials of Stochastic Finance*. World Scientific, River Edge, NJ. MR1695318
[3] SHIRYAEV, A. N., KABANOV, Y. M., KRAMKOV, D. O. and MELNIKOV, A. V. (1994). To the theory of computations of European and American options I. Discrete time. *Theory Probab. Appl.* **39** 14–60. MR1348190

INSTITUTE OF MATHEMATICS
THE HEBREW UNIVERSITY
JERUSALEM 91904
ISRAEL
E-MAIL: kifer@math.huji.ac.il